\documentclass{mcom-l}

\usepackage{amssymb,amsthm,amsfonts,amsmath,amsbsy,amstext,amsxtra,amsopn,graphicx,array,url,cite,latexsym}

\newcounter{lineno}
\newcounter{linelineno}[lineno]
\newcounter{linelinelineno}[linelineno]

\newcounter{algnum}
\renewcommand{\thealgnum}{\arabic{algnum}}

\newenvironment{algorithm}[3]
{
        \refstepcounter{algnum}
        \bigskip\goodbreak\hrule\medskip
                \leftline{{\sc Algorithm \thealgnum:} \bf #1}
                \medskip\hrule\medskip\nobreak
                \tt
        \begin{tabbing}
        OUTPUT: \= \kill
        INPUT:  \>#2 \\
        OUTPUT: \>#3
        \end{tabbing}
                \begin{tabbing}
10.\ \=word\=word\=word\=word\=word\=word\=\kill
                \smallskip \setcounter{lineno}{0}
}
{           \end{tabbing}
            \par\nobreak\hrule\bigskip
            \rm
}

\newcommand{\nline}{\refstepcounter{lineno} %
\'\thelineno.\>%
}

\newcommand{\nnline}{\refstepcounter{lineno}%
\'\thelineno.\>\>%
}

\newtheorem{theorem}{Theorem}[section]
\newtheorem{lemma}[theorem]{Lemma}
\newtheorem{heuristic}[theorem]{Heuristic}

\theoremstyle{definition}

\theoremstyle{remark}

\numberwithin{equation}{section}

\newcommand{\Z}{\mathbb{Z}}

\newcommand{\F}{\mathbb{F}}
\newcommand{\Q}{\mathbb{Q}}

\newcommand{\C}{\mathbb{C}}
\newcommand{\ie}{i.e., }

\newcommand{\ea}{{\em et al.}}

\newcommand{\ord}{\operatorname{ord}}

\begin{document}

\title[An efficient deterministic test for Kloosterman sum zeros]{An efficient deterministic test for Kloosterman sum zeros}

\author{Omran~Ahmadi}
\address{Claude Shannon Institute, UCD CASL, University College Dublin, Ireland}
\email{omran.ahmadi@ucd.ie}

\author{Robert~Granger}
\address{Claude Shannon Institute, Dublin City University, Ireland}
\curraddr{Claude Shannon Institute, UCD CASL, University College Dublin, Ireland}
\email{rgranger@computing.dcu.ie}
\thanks{Both authors are supported by the Claude Shannon Institute, Science Foundation Ireland Grant No. 06/MI/006.}


\begin{abstract}
We propose a simple deterministic test for deciding whether or not an
element $a \in \F_{2^n}^{\times}$ or $\F_{3^n}^{\times}$ is a zero of the corresponding
Kloosterman sum over these fields, and rigorously analyse its runtime. 
The test seems to have been overlooked in the literature.
The expected cost of the test for binary fields is a single
point-halving on an associated elliptic curve, while for ternary fields
the expected cost is one half of a point-thirding on an associated elliptic curve.
For binary fields of practical interest, this represents an $O(n)$ speedup over
the previous fastest test. By repeatedly invoking the test on random elements 
of $\F_{2^n}^{\times}$ we obtain the most efficient probabilistic
method to date to find non-trivial Kloosterman sum zeros. 
The analysis depends on the distribution of Sylow $p$-subgroups in the two
families of associated elliptic curves, which we ascertain using a theorem due to Howe.
\end{abstract}

\maketitle

\section{Introduction}
\label{sec:intro}

For a finite field $\F_{p^n}$, the Kloosterman sum $\mathcal{K}_{p^n}: \F_{p^n}
\rightarrow \C$ can be defined
by 
\[
\mathcal{K}_{p^n}(a) = 1 + \sum_{x \in \F_{p^n}^{\times}} \zeta^{\text{Tr}(x^{-1} + ax)},
\]
where $\zeta$ is a primitive $p$-th root of unity and Tr denotes the
absolute trace map $\text{Tr}:\F_{p^n} \rightarrow \F_p$, defined by 
\[
\text{Tr}(x) = x + x^p + x^{p^2} + \cdots + x^{p^{n-1}}. 
\]
Note that in some contexts the Kloosterman sum is defined to be just the
summation term without the added `$1$'~\cite{katz}. 
As one would expect, a Kloosterman (sum) zero is simply an element $a \in \F_{p^n}^{\times}$ for which
$\mathcal{K}_{p^n}(a) = 0$. 

Kloosterman sums have recently become the focus of much research, most notably
due to their applications in cryptography and coding theory
(see~\cite{gong,moisiocode} for example). 
In particular, zeros of $\mathcal{K}_{2^{n}}$ lead to bent
functions from $\F_{2^{2n}} \rightarrow \F_{2}$~\cite{dillon}, and similarly
zeros of $\mathcal{K}_{3^{n}}$ give rise to ternary bent functions~\cite{helleseth1}.

It was recently shown that zeros of Kloosterman sums only exist in
characteristics 2 and 3~\cite{kononen}, and hence these are the only cases we consider.
Finding such zeros is regarded as being difficult, and 
recent research has tended to focus on characterising Kloosterman sums modulo
small integers~\cite{charpin,moisio2,lisonek,lisonek2,lisonek3,faruk1,faruk2,faruk3,faruk4}. 
While these results are interesting in their own
right, they also provide a sieve which may be used to
eliminate elements of a certain form prior to testing whether they are
Kloosterman zeros or not, by some method.

It has long been known that Kloosterman sums over binary and ternary fields 
are intimately related to the group orders of members of two families of elliptic curves over
these fields~\cite{katz,wolfmann,moisio,geer}. In particular, for $p \in \{2,3\}$
the Kloosterman sum $\mathcal{K}_{p^n}(a)$ is equal to one minus the trace
of the Frobenius endomorphism of an associated elliptic curve $E_{p^n}(a)$. 
As such, one may use $p$-adic methods --- originally due to
Satoh~\cite{satoh} --- to compute the group orders of these elliptic curves,
and hence the corresponding Kloosterman sums.
The best $p$-adic point counting method asymptotically takes $O(n^2 \log^2{n}\log\log{n})$ bit
operations and requires $O(n^2)$ memory; see Vercauteren's thesis~\cite{frethesis} for
contributions and a comprehensive survey. 

Rather than count points, Lison\v{e}k has suggested that if instead one only wants to check 
whether a given element is a zero, one can do so by testing whether a random point of 
$E_{p^n}(a)$ has order $p^n$, via point multiplication~\cite{lisonek}.
Asymptotically, this has a similar bit complexity to
the point counting approach, requires less memory, but
is randomised. For fields of practical interest, it is reported that this
approach is superior to point counting~\cite[\S3]{lisonek}, and  
using this method Lison\v{e}k was able to find a
zero of $\mathcal{K}_{2^n}$ for $n \le 64$ and 
$\mathcal{K}_{3^n}$ for $n \le 34$, in a matter of days.

In this paper we take the elliptic curve connection to a logical
conclusion, in terms of proving divisibility results of Kloosterman sums by
powers of the characteristic. In particular we give an efficient deterministic algorithm to compute the 
Sylow $2$- and $3$-subgroups of the associated elliptic curves in
characteristics $2$ and $3$ respectively, along with a generator (these
subgroups are cyclic in the cases considered). Moreover, the average case runtimes of the
two algorithms are rigorously analysed. For binary fields of practical
interest, the test gives an $O(n)$ speedup over the point multiplication test.

Finding a single Kloosterman zero --- which is often all that
is needed in applications --- is then a matter of testing random field elements until
one is found, the success probability of which crucially depends on the number of Kloosterman zeros, see~\cite{katz}
and~\S\ref{exactformula}. Our runtime analysis provides a non-trivial
upper bound on this number, and consequently finding a Kloosterman zero with this
approach still requires time exponential in the size of the field.
We note that should one want to find {\em all}\hspace{0.5mm} Kloosterman
zeros over $\F_{2^n}$, rather than just one, then one can use the fast Walsh-Hadamard transform
(see~\cite{fft} for an overview), which requires $O(2^n \cdot n^2)$ bit operations and $O(2^n \cdot n)$ space.

The sequel is organised as follows. In \S\ref{connection} we detail the basic
connection between Kloosterman sums and two families of elliptic
curves. In \S\ref{determine} we present the main idea behind our algorithm, while 
\S\ref{binary} and \S\ref{ternary} explore its specialisation to binary and
ternary fields respectively. In \S\ref{noofits} we present data on the runtime
of the two algorithms, provide a heuristic analysis which attempts to explain the
data, and give an exact formula for the average case runtime. In \S\ref{mainresult}
we rigorously prove the expected runtime, while in \S\ref{compare} we assess
the practical efficiency of the tests. We finally make some concluding remarks 
in \S\ref{conc}.


\section{Connection with elliptic curves}\label{connection}

Our observations stem from the following three simple lemmas, which connect
Kloosterman sums over $\F_{2^n}$ and $\F_{3^n}$ with the group orders of 
elliptic curves in two corresponding families. The first is due to Lachaud and
Wolfmann~\cite{wolfmann}, the second Moisio~\cite{moisio}, while the third was 
proven by Lison\v{e}k~\cite{lisonek}.

\begin{lemma}\label{lis1}
Let $a \in \F_{2^n}^{\times}$ and define the elliptic curve $E_{2^n}(a)$ over
$\F_{2^n}$ by 
\[
E_{2^n}(a): y^2 + xy = x^3 + a.
\]
Then $\#E_{2^n}(a) = 2^n + \mathcal{K}_{2^n}(a)$.
\end{lemma}

\begin{lemma}\label{lis2}
Let $a \in \F_{3^n}^{\times}$ and define the elliptic curve $E_{3^n}(a)$ over
$\F_{3^n}$ by 
\[
E_{3^n}(a): y^2 = x^3 + x^2 - a.
\]
Then $\#E_{3^n}(a) = 3^n + \mathcal{K}_{3^n}(a)$.
\end{lemma}

\begin{lemma}\label{lis3}
Let $p \in \{2, 3\}$, let $a \in \F_{p^n}^{\times}$, and let $1 \leq h \leq
n$. Then $p^h \mid \mathcal{K}_{p^n}(a)$ if and only if there exists a point
of order $p^h$ on $E_{p^n}(a)$.
\end{lemma}

Lemma~\ref{lis3} is a simple consequence of the structure theorem for elliptic curves
over finite fields. Note that for $p \in \{2,3\}$, by Lemmas~\ref{lis1} and~\ref{lis2} we have 
$\mathcal{K}_{p^n}(a) = 0$ if and only if $E_{p^n}(a)$ has order $p^n$. By Lemma~\ref{lis3}, this is equivalent
to $E_{p^n}(a)$ having a point of order $p^n$, and hence finding a point of order $p^n$
proves that $\mathcal{K}_{p^n}(a) = 0$, since $p^n$ is the only element
divisible by $p^n$ in the Hasse interval. For the remainder of the
paper, when we refer to a prime $p$ we implicitly presume $p \in \{2,3\}$.


\section{Determining the Sylow $p$-subgroup of $E_{p^n}(a)$}\label{determine}

It is easy to show that $\mathcal{K}_{2^n}(a) \equiv 0 \pmod{4}$ and
$\mathcal{K}_{3^n}(a) \equiv 0 \pmod{3}$ for all $a \in \F_{2^n}^{\times}$ and $\F_{3^n}^{\times}$ respectively.
One way to see this is to observe that $E_{2^n}(a)$ possesses a point of order
$4$ (see~\S\ref{binary}) and $E_{3^n}(a)$ possesses a point of order $3$
(see~\S\ref{ternary}), and hence by Lagrange's theorem, $4 \mid \#E_{2^n}(a)$ and $3 \mid \#E_{3^n}(a)$. 

For an integer $x$, let $\ord_p(x)$ be the exponent of the maximum power of $p$ that divides $x$.
For $a \in \F_{p^n}^{\times}$, let $h = \ord_p(\#E_{p^n}(a))$. By
Lemma~\ref{lis3} the Sylow $p$-subgroup $S_{p}(E_{p^n}(a))$ is cyclic of order
$p^h$, and hence has $(p-1)p^{h-1}$ generators. Multiplying
these by $p$ results in the $(p-1)p^{h-2}$ generators of the order $p^{h-1}$
subgroup. Continuing this multiplication by $p$ process, after $h-1$ steps
one arrives at the $p$-torsion subgroup $E_{p^n}(a)[p]$, consisting of $p-1$ order-$p$ points
and the identity element $\mathcal{O}$. These considerations reveal the
structure of the $p$-power torsion subgroups $E_{p^n}(a)[p^k]$ for $1 \le k
\le h$, which one may view as a tree, with $\mathcal{O}$ as the root node. The root has $p-1$ children which are the
non-identity points in $E_{p^n}(a)[p]$. If $h>1$ each of these $p-1$ nodes 
has $p$ children: the elements of $E_{p^n}(a)[p^2] \setminus E_{p^n}(a)[p]$. 
For $1 < k < h$, each of the $(p-1)p^{k-1}$ depth-$k$ nodes have
$p$ children, while at depth $h$ we have $(p-1)p^{h-1}$ leaf nodes.

Using a division polynomial approach Lison\v{e}k was able to prove a necessary
condition on $a \in \F_{2^n}^{\times}$ such that $\mathcal{K}_{2^n}(a)$ is divisible by $16$, and likewise a
necessary condition on $a \in \F_{3^n}^{\times}$ such that  
$\mathcal{K}_{3^n}(a)$ is divisible by $9$. While necessary conditions for 
the divisibility of $\mathcal{K}_{2^n}(a)$ by $2^k$ have since been derived
for $k \le 8$~\cite{faruk3}, and for the divisibility of $\mathcal{K}_{3^n}(a)$
by $3^k$ for $k \le 3$~\cite{faruk4}, these use $p$-adic
methods; the division polynomial approach seemingly being too cumbersome to progress
any further.

However, the process outlined above --- taking a generator of $S_{p}(E_{p^n}(a))$
and multiplying by $p$ repeatedly until the non-identity elements of the
$p$-torsion are obtained --- can be reversed, easily and efficiently, using
point-halving in even characteristic, and point-thirding in characteristic
three, as we demonstrate in the ensuing two sections. Furthermore, due to
the cyclic structure of $S_{p}(E_{p^n}(a))$, at each depth, either all points are 
divisible by $p$, or none are. This means one can determine the height of the tree 
by using a depth-first search, without any backtracking; in particular, when a point $P$
at a given depth can not be halved or thirded, this depth is $\log_{p}(|S_{p}(E_{p^n}(a))|)$, and
$P$ is a generator. Furthermore, one can do this without
ever computing the group order of the curve. 

This process has been considered previously by Miret~\ea, for determining the Sylow $2$-subgroup of
elliptic curves over arbitrary finite fields of characteristic $ > 2$~\cite{miret1}; for $p=2$ the 
algorithm follows easily from the above considerations and point-halving, which is well studied in cryptographic
circles~\cite{knudsen,schroeppel,omran}, and is known to be more than
twice as fast as point-doubling in some cases~\cite{handbook}. 
For primes $l > 2$, Miret~\ea \ also addressed how to compute the Sylow $l$-subgroup of elliptic
curves over arbitrary finite fields provided that $l$ was not the
characteristic of the field~\cite{miret2}. Therefore we address here the case
$l = p = 3$, for the family of curves $E_{3^n}(a)$. 

We summarise this process in Algorithm~\ref{sylow}. Regarding notation, we say that a
point $P$ is $p$-divisible if there exists a point $Q$ such that $[p]Q = P$, and
write $Q = [1/p]P$.

\begin{algorithm}{DETERMINE $S_{p}(E_{p^n}(a))$}
{$a \in \F_{p^n}^{\times}$, $P \in E_{p^n}(a)[p] \setminus \{\mathcal{O}\}$}
{$(h,P_h)$ where $h = \ord_{p}(\#E_{p^n}(a))$ and $\langle P_h \rangle = S_{p}(E_{p^n}(a))$}
\label{sylow}
\nline $\text{counter} \leftarrow 1$; \\
\nline While $P$ is $p$-divisible do: \\
\nnline P := [1/p]P; \\ 
\nnline counter++; \\
\nline Return $(\text{counter},P)$
\end{algorithm}

Observe that Algorithm~\ref{sylow} is deterministic, provided that a
deterministic method of dividing a $p$-divisible point by $p$ is fixed once and for all, which we do
for $p=2$ and $p=3$ in~\S\ref{binary} and~\S\ref{ternary} respectively. For a given
field extension under consideration, choosing an appropriate field
representation and basis can also be performed deterministically, via
sequential search, however we consider this to be part of the setup phase
and do not incorporate setup costs when assessing the runtime of Algorithm~\ref{sylow}.


\section{Binary fields}\label{binary}

We now work out the details of Algorithm~\ref{sylow} for the family of curves
$E_{2^n}(a)$. For a fixed $n$, given a point $P =(x,y) \in E_{2^n}(a)$, $[2]P = (\xi,\eta)$
is given by the formula:
\begin{eqnarray}
\nonumber \lambda &=& x + y/x,\\
\label{half} \xi &=& \lambda^2 + \lambda,\\
\nonumber \eta &=& x^2 + \xi(\lambda+1).
\end{eqnarray}
To halve a point, one needs to reverse this process, \ie given $Q=(\xi,\eta)$, find (if
possible) a $P=(x,y) \in E_{2^n}(a)$ such that $[2]P = Q$. To do so, one first
needs to solve~(\ref{half}) for $\lambda$, which has a solution in
$\F_{2^n}$ if and only if $\text{Tr}(\xi) = 0$, since the trace of the
right-hand side is zero for every $\lambda \in \F_{2^n}$, and one can provide
an explicit solution in this case, as detailed in \S\ref{solvequad}. Observe that if
$\lambda$ is a solution to~(\ref{half}) then so is $\lambda + 1$. 
Assuming $\lambda$ has been computed, one then has 
\begin{eqnarray}
\nonumber x &=& (\eta + \xi(\lambda+1))^{1/2}, \\
\nonumber y &=& x(x+\lambda),
\end{eqnarray}
which for the two choices of $\lambda$ gives both points whose duplication is
$Q=(\xi,\eta)$. 

Aside from the cost of computing $\lambda$, the computation of $P = (x,y)$ as
above requires two field multiplications. As detailed in
Algorithm~\ref{sylow2}, this can be reduced to just one by
using the so-called $\lambda$-representation of a
point~\cite{knudsen,schroeppel}, where an affine point $Q = (\xi,\eta)$ is
instead represented by $(\xi,\lambda_{Q})$, with
\[
\lambda_Q = \xi + \frac{\eta}{\xi}.
\]
In affine coordinates, there is a unique $2$-torsion point $(0,a^{1/2})$, 
which halves to the two order $4$ points $P_{4}^{+} = (a^{1/4},a^{1/2})$,
$P_{4}^{-} = (a^{1/4},a^{1/2} + a^{1/4})$. The corresponding
$\lambda$-representations of each of these are $(a^{1/4},0)$ and $(a^{1/4},1)$ respectively.
For simplicity, we choose to use the former as the starting point in Algorithm~\ref{sylow2}.

\begin{algorithm}{DETERMINE $S_{2}(E_{2^n}(a))$}
{$a \in \F_{2^n}^{\times}$, $(x=a^{1/4}, \lambda = 0)$}
{$(h,P_h)$ where $h = \ord_{2}(\#E_{2^n}(a))$ and $\langle P_h \rangle = S_{2}(E_{2^n}(a))$}
\label{sylow2}
\nline $\text{counter} \leftarrow 2$; \\
\nline While $\text{Tr}(x) = 0$ do: \\
\nnline Solve $\widehat{\lambda}^2 + \widehat{\lambda} + x = 0$; \\
\nnline $t \leftarrow x(x + \lambda + \widehat{\lambda})$; \\
\nnline $x \leftarrow \sqrt{t}$; \\
\nnline $\lambda \leftarrow \widehat{\lambda} + 1$; \\
\nnline counter++; \\
\nline Return $(\text{counter},P = (x,x(x + \lambda)))$
\end{algorithm}

Observe that if the $x$-coordinate $a^{1/4}$ of $P_{4}^{\pm}$ satisfies
$\text{Tr}(a^{1/4}) = \text{Tr}(a) = 0$, then there exist four points of order
$8$, and hence $8 \mid \mathcal{K}_{2^n}(a)$, which was first observed by van
der Geer and van der Vlugt~\cite{geer}, and later by several others~\cite{helleseth2,charpin2,lisonek}.

\subsection{Solving $\widehat{\lambda}^2 + \widehat{\lambda} + x = 0$}\label{solvequad}

For odd $n$, let $\widehat{\lambda}$ be given by the following function,
which is known as the {\em half trace}:
\begin{equation}\label{halftrace}
\widehat{\lambda}(x) = \sum_{i=0}^{(n-1)/2} x^{2^{2i}}.
\end{equation}
One can easily verify that this $\widehat{\lambda}$ satisfies the stated
equation. When $n$ is even, the half trace approach will not
work, essentially because $\text{Tr}_{\F_{2^{n}}/\F_2}(1) = 0$.
Hence fix an element $\delta \in \F_{2^n}$ with $\text{Tr}_{\F_{2^{n}}/\F_2}(\delta) = 1$. Such a
$\delta$ can be found during the setup phase via the sequential search of the trace of
the polynomial basis elements, or by using the methods of~\cite{omran}. A solution to
equation~(\ref{half}) is then given by~\cite[Chapter II]{ECC1}:
\begin{equation}\label{fasthalf}
\widehat{\lambda}(x) = \sum_{i=0}^{n-2} \bigg( \sum_{j=i+1}^{n-1} \delta^{2^j} \bigg) x^{2^i},
\end{equation}
as may be verified. Note that for odd $n$, $\delta = 1$ suffices and so~(\ref{fasthalf})
simplifies to~(\ref{halftrace}). The inner sums of equation~(\ref{fasthalf}) can be precomputed, 
and for a general $\delta \in \F_{2^n}$ the computation of $\widehat{\lambda}(x)$ would require $n-1$ 
multiplications in $\F_{2^n}$, which together with the multiplication coming
from \url{line} \url{4} of Algorithm~\ref{sylow2}, gives a total of $n$ full $\F_{2^n}$-multiplications.

However, should $\F_{2^n}$ contain a subfield of odd index, then one
can reduce this cost as follows. Let $n = 2^m n'$ with
$m \ge 1$ and $n'$ odd. Constructing $\F_{2^n}$ as a degree $n'$ extension of
$\F_{2^{2^m}}$, fix a $\delta \in \F_{2^{2^m}}$ with
$\text{Tr}_{\F_{2^{2^m}}/\F_{2}}(\delta)=1$. Then 
\[
\text{Tr}_{\F_{2^{2^m\cdot n'}}/\F_2}(\delta)  = n' \cdot
\text{Tr}_{\F_{2^{2^m}}/\F_2}(\delta) = 1. 
\]
Hence this $\delta$ can be used in~(\ref{fasthalf}).
As $\delta^{2^{2^m}} = \delta$, upon
expanding~(\ref{fasthalf}) in terms of
$\{\delta^{2^0},\delta^{2^1},\ldots,\delta^{2^{2^m - 1}}\}$, we see that
at most $2^m$ multiplications of elements of $\F_{2^{2^m}}$ by elements of
$\F_{2^n}$ are required. So the smaller the largest power of $2$ dividing $n$ is, the
faster one can compute $\widehat{\lambda}(x)$.

However, since the expressions for $\widehat{\lambda}(x)$ in~(\ref{halftrace})
and~(\ref{fasthalf}) are linear maps, in practice it is far more efficient for
both odd and even $n$ to precompute and store
$\{\widehat{\lambda}(t^i)\}_{i=0,\ldots,n-1}$ during setup, 
where $\F_{2^n} = \F_{2}(t)$ and $x = \sum_{i=0}^{n-1} x_it^i$. One then has
\[
\widehat{\lambda}(x) = \sum_{i=0}^{n-1} x_i \widehat{\lambda}(t^i).
\]
On average just $n/2$ additions in $\F_{2^n}$ are required for each point-halving. Both the
storage required and execution time can be further reduced~\cite{handbook}. 
We defer consideration of the practical efficiency of
Algorithm~\ref{sylow2} until~\S\ref{compare2}.


\section{Ternary fields}\label{ternary}

Let $Q=(\xi,\eta) \in E_{3^n}(a)$. To find $P=(x,y)$ such that $[3]P = Q$, when possible, we do the following.
As in~\cite[\S4]{miret2}, we have
\[
x([3]P) = x(P) - \frac{\Psi_{2}(x,y)\Psi_{4}(x,y)}{\Psi_{3}^{2}(x,y)},
\]
or
\[
(x - \xi)\Psi_{3}^{2}(x,y) - \Psi_{2}(x,y)\Psi_{4}(x,y) = 0,
\]
where $\Psi_l$ is the $l$-th division polynomial.
Working modulo the equation of $E_{3^n}(a)$, this becomes
\[
x^9 - \xi x^6 + a(1-\xi)x^3 - a^2(a + \xi) = 0,
\]
whereupon substituting $X = x^3$ gives
\begin{equation}\label{3div}
f(X) = X^3 - \xi X^2 + a(1-\xi)X - a^2(a + \xi) = 0.
\end{equation}
To solve~(\ref{3div}), we make the transformation
\[
g(X) = X^3 f\bigg(\frac{1}{X} - \frac{a(1-\xi)}{\xi}\bigg) = \frac{a^2 \eta^2}{\xi^3} X^3 -\xi X  + 1. 
\]
Hence we must solve
\[
X^3 - \frac{\xi^4}{a^2 \eta^2} X  + \frac{\xi^3}{a^2 \eta^2} = 0.
\]
Writing $X = \frac{\xi^2}{a \eta}\widehat{X}$ this becomes
\begin{equation}\label{3tri}
\widehat{X}^3 - \widehat{X} + \frac{a \eta}{\xi^3} = 0.
\end{equation}
Our thirding condition is then simply $\text{Tr}(a \eta/\xi^3) = 0$,
since as in the binary case, for every element $\widehat{X} \in \F_{3^n}$ we have $\text{Tr}(\widehat{X}^3
- \widehat{X}) = 0$, and if so then one can provide an explicit
solution, as detailed in \S\ref{solvecube}. Observe that if $\widehat{X}$ is a solution to~(\ref{3tri})
then so is $\widehat{X} \pm 1$.
Unrolling the transformations leads to the following algorithm, with input the
$3$-torsion point $P_3 = (a^{1/3},a^{1/3})$.

\begin{algorithm}{DETERMINE $S_{3}(E_{3^n}(a))$}
{$a \in \F_{3^n}^{\times}$, $(x=a^{1/3}, y = a^{1/3})$}
{$(h,P_h)$ where $h = \ord_3(\#E_{3^n}(a))$ and $\langle P_h \rangle = S_{3}(E_{3^n}(a))$}
\label{sylow3}
\nline $\text{counter} \leftarrow 1$; \\
\nline While $\text{Tr}(ay/x^3) = 0$ do: \\
\nnline Solve $\widehat{X}^3 - \widehat{X} + \frac{ay}{x^3} = 0$; \\
\nnline $x \leftarrow \bigg(\frac{ay}{x^2\widehat{X}} - \frac{a(1-x)}{x} \bigg)^{1/3}$; \\
\nnline $y \leftarrow \big(x^3 + x^2 - a\big)^{1/2}$; \\
\nnline counter++; \\
\nline Return $(\text{counter},P = (x,y))$
\end{algorithm}

Observe that as with Algorithm~\ref{sylow2}, if the point $P_3$ satisfies
$\text{Tr}(a \cdot a^{1/3}/a) = \text{Tr}(a) = 0$, then there is a point of order
$9$, and hence $9 \mid \mathcal{K}_{3^n}(a)$, which again was first proven
in~\cite{geer}, and later by others~\cite{lisonek,faruk1}.

\subsection{Solving $\widehat{X}^3 - \widehat{X} + \frac{a y}{x^3} = 0$}\label{solvecube}

Let $\beta = \frac{a y}{x^3}$, and let $\delta \in \F_{3^n}$ be an
element with $\text{Tr}_{\F_{3^{n}}/\F_3}(\delta) = 1$, which can be found deterministically
during the setup phase. It is then a simple matter to verify that 
\begin{equation}\label{fastthird}
\widehat{X}(\beta) = \sum_{i=0}^{n-2} \bigg( \sum_{j=i+1}^{n-1} \delta^{3^j} \bigg) \beta^{3^i}
\end{equation}
is a solution to equation~(\ref{3tri}).

For $n \equiv 1 \pmod{3}$, one may choose $\delta = 1$ and the expression for $\widehat{X}(\beta)$ in
equation~(\ref{fastthird}) simplifies to 
\[
\widehat{X}(\beta) = \sum_{i=1}^{(n-1)/3} \left(\beta^{3^{3i-1}} - \beta^{3^{3i-2}}\right). 
\]

For $n \equiv 2 \pmod{3}$, one may choose $\delta = -1$ and the expression for $\widehat{X}(\beta)$ in
equation~(\ref{fastthird}) simplifies to 
\[
\widehat{X}(\beta) = -\beta + \sum_{i=1}^{(n-2)/3} \left(\beta^{3^{3i-1}} - \beta^{3^{3i}}\right). 
\]

For $n \equiv 0 \pmod{3}$, one can use the approach described in
\S\ref{solvequad} to pick $\delta$ from the smallest subfield of $\F_{3^n}$ of index
coprime to $3$, in order to reduce the cost and the number of 
multiplications required to solve~(\ref{3tri}). As in the binary
case, one can also exploit the linearity of $\widehat{X}(\beta)$ and
precompute and store $\{\widehat{X}(t^i)\}_{i=0,\ldots,n-1}$ during setup,
where $\F_{3^n} = \F_{3}(t)$ and $\beta = \sum_{i=0}^{n-1} \beta_it^i$, in order to reduce the cost of
solving~(\ref{3tri}) to an average of $2n/3$ additions.
We defer consideration of the practical efficiency of Algorithm~\ref{sylow3} until~\S\ref{compare3}.


\section{Heuristic analysis of the expected number of iterations}\label{noofits}

For any input $a \in \F_{p^n}^{\times}$, the runtime of Algorithm~\ref{sylow} is proportional to 
the number of loop iterations performed, which is precisely the height of the corresponding
Sylow $p$-subgroup tree, $h = \log_{p}(|S_p(E_{p^n}(a))|)$.
In this section we present experimental data for the distribution of these
heights for $p \in \{2,3\}$, provide a heuristic argument to explain them,
and give an exact formula for the average case runtime.
Since we are interested in the average number of loop iterations\footnote{The worst case
being $n$ iterations, which of course is the best case when searching for a
Kloosterman zero.}, we consider the arithmetic mean of the heights of the 
Sylow $p$-subgroup trees, or equivalently the logarithm of the geometric mean of their
orders.

\subsection{Experimental data}

In order to gain an idea of how $\{\log_{p}(|S_p(E_{p^n}(a))|)\}_{a \in  \F_{p^n}^{\times}}$ is distributed, 
we computed all of them for several small extensions of $\F_p$. Tables~\ref{dist2} and~\ref{dist3}
give the results for $p=2$ and $p=3$ respectively. 

Observe that for $p=2$, the first two columns are simply $2^n - 1 =
|\F_{2^n}^{\times}|$, reflecting the fact that all of the curves
$\{E_{2^n}(a)\}_{a \in \F_{2^n}^{\times}}$ have order divisible by
$4$. Similarly for $p=3$, the first column is given by $3^n-1 = |\F_{3^n}^{\times}|$,
reflecting the fact that all the curves $\{E_{3^n}(a)\}_{a \in \F_{3^n}^{\times}}$ have order divisible by
$3$. Furthermore, since exactly half of the elements of $\F_{2^n}$ have zero trace,
the third column for $p=2$ is given by $2^{n-1}-1$. Likewise for $p=3$, the second column
is given by $3^{n-1} - 1$, since exactly one third of the elements of $\F_{3^n}$ have zero trace.
For $p=2$ there is an elegant result due Lison\v{e}k and Moisio which gives a
closed formula for the $n$-th entry of column $4$ of
Table~\ref{dist2}~\cite[Theorem 3.6]{lisonek3}, which includes the $a=0$ case,
namely:
\begin{equation}\label{column4}
(2^n - (-1 + i)^n - (-1 - i)^n)/4. 
\end{equation}
Beyond these already-explained columns, it appears that as one successively moves one column
to the right, the number of such $a$ decreases by an approximate factor of $2$ or $3$
respectively, until the number of Kloosterman zeros is reached, which by
Hasse bound occurs as soon as  $p^k > 1 + 2p^{n/2}$, or $k > n/2 + \log_{p}2$.

\begin{table}\caption{$\# \{E_{2^n}(a)\}_{a \in \F_{2^n}^{\times}}$ whose group order is divisible by $2^{k}$}\label{dist2}
\begin{center}
\begin{tabular}{|c|r|r|r|r|r|r|r|r|r|r|r|r|r|}
\hline
$n \backslash k$ & $1$ & $2$    & $3$   & $4$   & $5$  & $6$  & $7$ & $8$ &
$9$ & $10$ & $11$ & $12$ & $13$\\
\hline
1  &  1       &   1     &        &       &       &      &      &     &  & & & &    \\
2  &  3       &   3    &        &       &       &      &      &     &   & & & &   \\
3  &  7       &   7    &     3  &       &       &      &      &     &  & & & &    \\
4  &  15      &   15   &     7  &   5   &       &      &      &     &  & & & &    \\
5  &  31      &   31   &    15  &   5   &   5   &      &      &     &  & & & &    \\
6  &  63      &   63   &    31  &  15   &  12   &  12  &      &     &   & & & &   \\
7  & 127      &  127   &    63  &  35   &  14   &  14  & 14   &     &   & & & &   \\
8  & 255      &  255   &   127  &  55   &  21   &  16  & 16   & 16  &  & & & &    \\
9  & 511      &  511   &   255  & 135   &  63   &  18  & 18   & 18  & 18 & & & &\\
10 & 1023     & 1023   &   511  & 255   & 125   &  65  & 60   & 60  & 60 & 60 & & &\\
11 & 2047     & 2047   &  1023  & 495   & 253   & 132  & 55   & 55  & 55 & 55 & 55 & & \\
12 & 4095     & 4095   &  2047  & 1055  & 495   & 252  & 84   & 72  & 72 & 72 & 72 & 72 &\\
13 & 8191     & 8191   &  4095  & 2015  & 1027  & 481  & 247  & 52  & 52 & 52
& 52 & 52 &  52\\
\hline
\end{tabular}
\end{center}
\end{table}

\begin{table}\caption{$\# \{E_{3^n}(a)\}_{a \in \F_{3^n}^{\times}}$ whose group order is divisible by $3^k$}\label{dist3}
\begin{center}
\begin{tabular}{|c|r|r|r|r|r|r|r|r|r|r|r|}
\hline
$n \backslash k$ & $1$ & $2$    & $3$   & $4$   & $5$  & $6$  & $7$ & $8$ &
$9$ & $10$ & $11$\\
\hline
1   & 2       &        &        &       &       &      &      &     &  & &   \\
2   & 8       & 2      &        &       &       &      &      &     &   & &   \\
3   & 26      & 8      & 3      &       &       &      &      &     &   & &   \\
4   & 80      & 26     & 4      & 4     &       &      &      &     &    & &  \\
5   & 242     & 80     & 35     & 15    & 15    &      &      &     &    & &  \\
6   & 728     & 242    & 83     & 24    & 24    & 24   &      &     &   & &   \\
7   & 2186    & 728    & 266    & 77    & 21    & 21   & 21   &     &   & &   \\
8   & 6560    & 2186   & 692    & 252   & 48    & 48   & 48   & 48  &   & &   \\
9   & 19682   & 6560   & 2168   & 741   & 270   & 108  & 108  & 108 & 108 & & \\
10  & 59048   & 19682  & 6605   & 2065  & 575   & 100  & 100  & 100 & 100 & 100 & \\
11  & 177146  & 59048  & 19547  & 6369  & 2596  & 924  & 264  & 264 & 264 & 264 & 264\\
\hline
\end{tabular}
\end{center}
\end{table}

\subsection{A heuristic for the expected number of iterations}\label{heuristic}

To explain the data in Tables~\ref{dist2} and~\ref{dist3}, we propose
the following simple heuristic (and prove the validity of its consequences in~\S\ref{mainresult}):

\begin{heuristic}\label{heur}
Over all $a \in \F_{p^n}^{\times}$, on any occurrence of {\em \url{line}} {\em \url{2}} 
of the loop in Algorithms~\ref{sylow2} and~\ref{sylow3}, regardless of the height of
the tree at that point, the argument of the $\F_{p^n}$ trace is uniformly
distributed over $\F_{p^n}$, and hence is zero with probability $1/p$. 
\end{heuristic}
While this assumption is clearly false at depths $> n/2 + \log_{p}2$, the data
in Tables~\ref{dist2} and~\ref{dist3} does support it (up to relatively small error terms).
In order to calculate the expected value of $\log_{p}(|S_p(E_{p^n}(a))|)$, 
we think of Algorithms~\ref{sylow2} and~\ref{sylow3} as running on all $p^n-1$ elements of
$\F_{p^n}^{\times}$ in parallel; we then sum the number of elements which
survive the first loop, then the second loop and the third loop etc., and divide this sum 
by $p^n-1$ to give the average. 
We now explore the consequences of Heuristic~\ref{heur}, treating the two
characteristics in turn.

For Algorithm~\ref{sylow2}, on the first occurrence of \url{line} \url{2},
$2^{n-1} - 1$ elements of $\F_{2^n}^{\times}$ have zero trace and hence
$2^{n-1}-1$ elements require an initial loop iteration.
On the second occurrence of \url{line} \url{2}, by Heuristic~\ref{heur},
approximately $2^{n-1}/2 = 2^{n-2}$ of the inputs have zero trace and so this number of
loop iterations are required. Continuing in this manner and summing over all
loop iterations at each depth, one obtains a total of
\[
2^{n-1} + 2^{n-1} + \cdots + 2 + 1 \approx 2^n,
\]
for the number of iterations that need to be performed for all $a \in
\F_{2^n}^{\times}$. Thus on average this is approximately one loop iteration
per initial element $a$. Incorporating the divisibility by $4$ of all curve orders, 
the expected value as $n \rightarrow \infty$ of $\log_{2}(|S_2(E_{2^n}(a))|)$ 
is $3$, and hence the geometric mean of $\{|S_2(E_{2^n}(a))|\}_{a \in
  \F_{2^n}^{\times}}$ as $n \rightarrow \infty$ is $2^3 = 8$.

For Algorithm~\ref{sylow3}, applying Heuristic~\ref{heur} and the same
reasoning as before, the total number of iterations required for all $a \in
\F_{3^n}^{\times}$ is
\[
3^{n-1} + 3^{n-2} + \cdots + 3 + 1 \approx 3^n/2.
\]
Thus on average this is approximately $1/2$ an iteration per initial element $a$, and
incorporating the divisibility by $3$ of all curve orders, the expected value as $n \rightarrow \infty$ 
of $\log_{3}(|S_3(E_{3^n}(a))|)$ is $3/2$, and hence the geometric 
mean of $\{|S_3(E_{3^n}(a))|\}_{a \in \F_{3^n}^{\times}}$ as $n \rightarrow \infty$ is $3^{3/2} = 3\sqrt{3}$.

\subsection{Exact formula for the average height of Sylow $p$-subgroup trees}\label{exactformula}

Let $p^n + t$ be an integer in the Hasse interval $I_{p^n} = [p^n + 1 - 2p^{n/2},p^n + 1 + 2p^{n/2}]$,
 which is assumed to be divisible by $4$ if $p=2$ and divisible by $3$ if $p=3$. 
Let $N(t)$ be the number of solutions in $\F_{p^n}^{\times}$ to $\mathcal{K}_{p^n}(a) = t$. 
The sum of the heights of the Sylow $p$-subgroup trees, over all $a \in \F_{p^n}^{\times}$, is
\begin{equation}\label{exactav}
T_{p^n} = \sum_{(p^n+t) \in I_{p^n}} N(t) \cdot \text{ord}_p(p^n+t),
\end{equation}
and thus the expected value of $\log_{p}(|S_p(E_{p^n}(a))|)$ is $T_{p^n}/(p^n-1)$.
The crucial function $N(t)$ in~(\ref{exactav}) has been evaluated by Katz and
Livn\'{e} in terms of class numbers~\cite{katz}. 
In particular, let $\alpha  = (t-1 + \sqrt{(t-1)^2 - 4p^n})/2$ for $t$ as above. Then 
\[
N(t) = \sum_{\text{orders} \ \mathcal{O}} h(\mathcal{O}), 
\]
where the sum is over all orders $\mathcal{O} \subset \Q(\alpha)$ which contain $\Z[\alpha]$.
It seems difficult to prove Heuristic~\ref{heur} or our implied estimates for $T_{p^n}$ using the Katz-Livn\'{e}
result directly.
However, using a natural decomposition of $T_{p^n}$ and a theorem due to
Howe~\cite{howe}, in the following section we show that the consequences of Heuristic~\ref{heur} as
derived in~\S\ref{heuristic} are correct.

\section{Main result}\label{mainresult}

We now present and prove our main result, which states that the expected value of
$\{\log_{p}(|S_p(E_{p^n}(a))|)\}_{a \in  \F_{p^n}^{\times}}$ is precisely as we derived heuristically in~\S\ref{heuristic}. 
To facilitate our analysis, for $1 \le k \le n$, we partition $T_{p^n}$ into the counting functions 
\begin{equation}\label{igusaprimer}
T_{p^n}(k) = \sum_{(p^n + t) \in I_{p^n}, p^k|(p^n + t)} N(t),
\end{equation}
so that by~(\ref{exactav}) we have
\begin{equation}\label{partition}
T_{p^n} = \sum_{k = 1}^{n} T_{p^n}(k). 
\end{equation} 
Indeed, the integers $T_{p^n}(k)$ are simply the $(n,k)$-th entries of
Tables~\ref{dist2} and~\ref{dist3} for $p=2$ and $3$ respectively, and thus
$T_{p^n}$ is the sum of the $n$-th row terms. Hence we already have $T_{2^n}(1) =
T_{2^n}(2) = 2^n-1$, $T_{2^n}(3) = 2^{n-1}-1$ and $T_{2^n}(4) = (2^n - (-1 +
i)^n - (-1 - i)^n)/4$ by~(\ref{column4}), and similarly $T_{3^n}(1) =
3^n -1$ and $T_{3^n}(2) = 3^{n-1}-1$.

\subsection{Estimating $T_{p^n}(k)$}

For $k \ge 2$, let $\mathcal{T}_{2^n}(k)$ be the set of $\F_{2^n}$-isomorphism classes of elliptic curves
$E/\F_{2^n}$ such that $\#E(\F_{2^n}) \equiv 0 \pmod{2^k}$. Similarly for $k
\ge 1$, let $\mathcal{T}_{3^n}(k)$ be the set of $\F_{3^n}$-isomorphism classes of elliptic curves
$E/\F_{3^n}$ such that $\#E(\F_{3^n}) \equiv 0 \pmod{3^k}$.
Observe that the elliptic curves $E_{2^n}(a)$ and $E_{3^n}(a)$ both have $j$-invariant
$1/a$~\cite[Appendix A]{Silverman}, and hence cover all the
$\overline{\F}_{2^n}$- and $\overline{\F}_{3^n}$-isomorphism
classes of elliptic curves over $\F_{2^n}$ and $\F_{3^n}$ respectively, except
for $j=0$. We have the following lemma.

\begin{lemma}\label{wouter}\cite[Lemma 6]{castryck}
Let $E/\F_q$  be an elliptic curve and let $[E]_{\F_q}$ be the set of $\F_q$-isomorphism
classes of elliptic curves that are $\overline{\F}_q$-isomorphic to $E$. Then for
$j \ne 0,1728$ we have $\#[E]_{\F_q} = 2$, and $[E]_{\F_q}$ consists of the
$\F_q$-isomorphism class of $E$ and the $\F_q$-isomorphism class of its quadratic twist $E^t$.
\end{lemma}

Let $\#E_{2^n}(a) = 2^n + 1 - t_a$, with $t_a$ the trace of Frobenius. Since
$j \ne 0$, by Lemma~\ref{wouter} the only other $\F_{2^n}$-isomorphism class with $j$-invariant $1/a$
is that of the quadratic twist $E_{2^n}^{t}(a)$, which has order $2^n + 1
+ t_a$. Since $t_a \equiv 1 \pmod{4}$, we have $\#E_{2^n}^{t}(a) \equiv 2
\pmod{4}$ and hence none of the $\F_{2^n}$-isomorphism classes of the quadratic twists of 
$E_{2^n}(a)$ for $a \in \F_{2^n}^{\times}$ are in $\mathcal{T}_{2^n}(k)$, for $k \ge 2$.
By an analogous argument, only the $\F_{3^n}$-isomorphism classes of $E_{3^n}(a)$ for $a \in \F_{3^n}^{\times}$
are in $\mathcal{T}_{3^n}(k)$, for $k \ge 1$. Furthermore, all curves
$E/\F_{2^n}$ and $E/\F_{3^n}$ with $j=0$ are
supersingular~\cite[\S3.1]{washington}, and therefore have group orders $\equiv 1 \pmod{4}$ and $\equiv
1 \pmod{3}$ respectively. Hence no $\F_{p^n}$-isomorphism classes of curves with
$j=0$ are in $\mathcal{T}_{p^n}(k)$ for $p \in \{2,3\}$.
As a result, for $2 \le k \le n$ we have 
\begin{equation}\label{Tequal}
|\mathcal{T}_{2^n}(k)| = T_{2^n}(k),
\end{equation}
and similarly, for $1 \le k \le n$ we have 
\begin{equation}
\nonumber |\mathcal{T}_{3^n}(k)| = T_{3^n}(k).
\end{equation}
Therefore in both cases, a good estimate for $|\mathcal{T}_{p^n}(k)|$ is all we need to 
estimate $T_{p^n}(k)$. The cardinality of $\mathcal{T}_{p^n}(k)$
is naturally related to the study of modular curves; in particular, considering
the number of $\F_{p^n}$-rational points on the Igusa curve of level $p^k$ allows one to prove 
Theorem~\ref{maintheorem} below~\cite{Igusa,Amilcar}. However, for simplicity (and generality) we
use a result due to Howe on the group orders of elliptic curves over finite fields~\cite{howe}. 
Consider the set 
\[
V(\F_q;N) = \{ E/\F_q: N \mid \#E(\F_q)\}\big/\cong_{\F_q}
\]
of equivalence classes of $\F_q$-isomorphic curves whose group orders are divisible by $N$. 
Following Lenstra~\cite{lenstra}, rather than estimate $V(\F_q;N)$ directly, Howe considers 
the weighted cardinality of $V(\F_q;N)$, where for a set $S$ of $\F_q$-isomorphism classes of elliptic
curves over $\F_q$, this is defined to be:
\[
\#'S = \sum_{[E] \in S} \frac{1}{\#\text{Aut}_{\F_q}(E)}.
\]
For $j \ne 0$ we have $\#\text{Aut}_{\overline{\F}_q}(E) =
2$~\cite[\S{III.10}]{Silverman} 
and since $\{\pm 1\} \subset \text{Aut}_{\F_q}(E)$ we have 
$\#\text{Aut}_{\F_q}(E) = 2$ also. Therefore, by the above discussion, 
for $p=2, k \ge 2$ and $p =3, k \ge 1$ we have
\begin{equation}\label{estimate}
|\mathcal{T}_{p^n}(k)| = 2 \cdot \#'V(\F_{p^n};p^k), 
\end{equation}
We now present Howe's result.
\begin{theorem}\label{howesthm}\cite[Theorem 1.1]{howe}
There is a constant $C \le  1/12 + 5\sqrt{2}/6 \approx 1.262$ such that the following
statement is true: Given a prime power $q$, let $r$ be the multiplicative arithmetic function
such that for all primes $l$ and positive integers $a$
\[
r(l^a) =
\begin{cases} 
\dfrac{1}{l^{a-1}(l-1)}, & \mbox{if } q \not\equiv 1
  \pmod{l^c};\\
\\
\dfrac{l^{b+1}+ l^b-1}{l^{a+b-1}(l^2-1)}, & \mbox{if } q \equiv 1
  \pmod{l^c},
\end{cases}
\]
where $b = \lfloor a/2 \rfloor$ and $c = \lceil a/2 \rceil$. Then for all
positive integers $N$ one has
\begin{equation}\label{howeformula}
\bigg| \frac{\#'V(\F_q;N)}{q} - r(N) \bigg| \le \frac{C N \rho(N)2^{\nu(N)}}{\sqrt{q}}, 
\end{equation}
where $\rho(N) = \prod_{p \mid N}((p+1)/(p-1))$ and $\nu(N)$ denotes the number
of distinct prime divisors of $N$.
\end{theorem}

Equipped with Theorem~\ref{howesthm}, we now present and prove our main theorem.

\begin{theorem}\label{maintheorem}
Let $p \in \{2,3\}$ and let $T_{p^n}(k)$ be defined as above. Then
\begin{itemize}
\item[(i)] For $3 \le k < n/4$ we have $T_{2^n}(k) = 2^{n - k + 2} + O(2^{k+n/2})$,
\item[(ii)] For $2 \le k < n/4$ we have $T_{3^n}(k) = 3^{n - k + 1} + O(3^{k+n/2})$,
\item[(iii)] $T_{2^n} = 3 \cdot 2^n + O(n \cdot 2^{3n/4})$,
\item[(iv)] $T_{3^n} = 3^{n+1}/2 + O(n \cdot 3^{3n/4})$,
\item[(v)] $\lim_{n \to \infty} T_{p^n}/(p^n-1) = 
\begin{cases}
3 \hspace{8mm} \text{if} \ p = 2,\\
3/2 \hspace{4.5mm} \text{if} \ p=3.
\end{cases}$
\end{itemize}
Furthermore, in $(i)-(iv)$ the implied constants in the $O$-notation are absolute and
effectively computable.
\end{theorem}

\begin{proof}
By equations~(\ref{Tequal}) and~(\ref{estimate}), and Theorem~\ref{howesthm}
with $l=p$, for $3 \le k \le n$ we have 
\[
\bigg| \frac{T_{2^n}(k)}{2^{n+1}} - \frac{1}{2^{k-1}} \bigg| \le \frac{C \cdot
  2^k\cdot 3 \cdot 2}{2^{n/2}},
\]
from which (i) follows immediately. Similarly for $2 \le k \le n$ we have 
\[
\bigg| \frac{T_{3^n}(k)}{2 \cdot 3^{n}} - \frac{1}{3^{k-1}\cdot 2} \bigg| \le \frac{C \cdot
  3^k \cdot (4/2) \cdot 2}{3^{n/2}},
\]
from which (ii) follows.
For (iii) we write equation~(\ref{partition}) as follows:
\begin{equation}
\nonumber T_{2^n} = \sum_{k = 1}^{n} T_{2^n}(k) = \sum_{k = 1}^{\lfloor n/4 \rfloor
  -1} T_{2^n}(k) + \sum_{k =\lfloor n/4 \rfloor}^{n} T_{2^n}(k).
\end{equation}
Freely applying (i), the first of the these two sums equals
\begin{eqnarray}
\nonumber & & 2^n + (2^n + 2^{n-1} + \cdots + 2^{n - \lfloor n/4 \rfloor + 2}) + O(2^{n/2 + 2}
+ 2^{n/2 + 3} + \cdots + 2^{n/2 + \lfloor n/4 \rfloor})\\
\nonumber &=& 2^n + 2^{n+1} - 2^{n - \lfloor n/4 \rfloor + 2} + O(2^{n/2 + \lfloor n/4 \rfloor + 1})\\
\nonumber &=& 2^n + 2^{n+1} + O(2^{3n/4}) = 3\cdot 2^n +  O(2^{3n/4}).
\end{eqnarray}
For the second sum, observe that $p^{k+1} \mid t \Longrightarrow p^k \mid t$ and so $T_{2^n}(k+1) \le T_{2^n}(k)$, which gives
\[
\sum_{k =\lfloor n/4 \rfloor}^{n} T_{2^n}(k) \le (3n/4 + 2) \cdot T_{2^n}(\lfloor
  n/4 \rfloor) = O(n \cdot 2^{3n/4}).
\]
Combining these two sums one obtains (iii). Part (iv) follows {\em mutatis mutandis}, which together with (iii) proves (v).
\end{proof}

Theorem~\ref{maintheorem} proves that for $k < n/4$, the distribution of the height function
$\log_{p}(|S_p(E_{p^n}(a))|)$ over $a \in \F_{p^n}^{\times}$ is approximately
geometric. Hence using an argument similar to the above one can prove that
asymptotically, the variance is $2$ for $p=2$, and $3/4$ for $p=3$.
Our proof also gives an upper bound on the number of Kloosterman zeros.
In particular, parts (i) and (ii) imply that for $k < n/4$, for increasing
$k$, $T_{p^n}(k)$ is decreasing, and hence the number of Kloosterman zeros is $O(p^{3n/4})$.
Shparlinski has remarked~\cite{shpar} that this upper bound follows from 
a result of Niederreiter~\cite{nied}, which refines an earlier result due to Katz~\cite{katzM}.
The Weil bound intrinsic to Howe's estimate fails to give any tighter bounds on $|T_{p^n}(k)|$ for $n/4  \le k \le
n/2$. Finding improved bounds on $|T_{p^n}(k)|$ for $k$ in this interval is an
interesting problem, since they would immediately give a better upper bound on the number of Kloosterman zeros.

While our proof only required the $l=p$ part of Howe's result
(when we could have used tighter bounds arising from an Igusa curve argument), 
the more general form, when combined with our approach, allows one to compute
the expected height of the Sylow $l$-subgroup trees for $l \ne p$ as well,
should this be of interest.

\section{Test Efficiency}\label{compare}

We now address the expected efficiency of Algorithms~\ref{sylow2} 
and~\ref{sylow3} when applied to random elements of $\F_{2^n}^{\times}$ and $\F_{3^n}^{\times}$ respectively.
Since the number of Kloosterman zeros is $O(p^{3n/4})$, by choosing random $a \in
\F_{p^n}^{\times}$ and applying our algorithms, one only has an exponentially small
probability of finding a zero. Hence we focus on those $n$ for which such
computations are currently practical and do not consider the asymptotic
complexity of operations.
For comparative purposes we first recall Lison\v{e}k's randomised Kloosterman
zero test~\cite{lisonek}.

\subsection{Lison\v{e}k's Kloosterman zero test}
For a given $a \in \F_{p^n}^{\times}$, Lison\v{e}k's test consists of taking 
a random point $P \in E_{p^n}(a)$, and computing $[p^n]P$ to see if it is
the identity element $\mathcal{O} \in E_{p^n}(a)$. If it is not, then by Lemmas~\ref{lis1} 
and~\ref{lis2} one has certified that the group order is not $p^n$ and thus $a$ is
not a Kloosterman zero. If $[p^n]P = \mathcal{O}$ and $[p^{n-1}]P \neq
\mathcal{O}$, then $\langle P \rangle = E_{p^n}(a)$ and $a$ is a Kloosterman
zero. In this case the probability that a randomly chosen point on the curve
is a generator is $1/2$ and $2/3$ for $p=2$ and $p=3$ respectively.
The test thus requires $O(n)$ point-doublings/triplings in $E_{2^n}(a)$ and
$E_{3^n}(a)$ respectively.

 
\subsection{Algorithm~\ref{sylow2} for $E_{2^n}(a)$}\label{compare2}

By Theorem~\ref{maintheorem}(v), only one loop iteration of
Algorithm~\ref{sylow2} is required on average. Each such iteration requires computing: a trace; solving~(\ref{half}); 
a multiplication; a square root; two additions; and a bit-flip.
This process has been extensively studied and optimised for point-halving in characteristic
$2$~\cite{handbook}. In particular, for $n=163$ and $n=233$, point-halving is
reported to be over twice as fast as point-doubling~\cite[Table 3]{handbook}.
Hence in this range of $n$, with a state-of-the-art implementation,
Algorithm~\ref{sylow2} is expected to be $\approx 2n$ times faster than
Lison\v{e}k's algorithm (or $\approx n$ times faster if for the latter one checks whether or not
$\text{Tr}(a) = 0$ before initiating the point multiplication).

For the field $\F_{2^{75}} = \F_2[t]/(t^{75} + t^6 + t^3 + t + 1)$, using a
basic MAGMA V2.16-12~\cite{magma} implementation of Algorithm~\ref{sylow2}, we found
the Kloosterman zero:
\begin{eqnarray*}
a & = & t^{74} + t^{73} + t^{68} + t^{67} + t^{66} + t^{65} + t^{63} + t^{62} + t^{59} + t^{57} + t^{56} + t^{55} + t^{52}\\ 
  & + & t^{44} + t^{43} + t^{41} + t^{40} + t^{39} + t^{38} + t^{37} + t^{36} + t^{35} + t^{34} + t^{31} + t^{30} + t^{29}\\
  & + & t^{28} + t^{25} + t^{24} + t^{23} + t^{22} + t^{19} + t^{16} + t^{15} + t^{14} + t^{13} + t^{12} + t^{11} + t^{8}\\ 
  & + & t^{7} + t^{6} + t^5 + t^4 + t^3 + t^2 + t,
\end{eqnarray*}
in 18 hours using eight AMD Opteron 6128 processors each running at 2.0
GHz. Due to MAGMA being general-purpose, without a built-in function for
point-halving, the above implementation has comparable efficiency to a full point multiplication by $2^{75}$ on $E_{p^n}(a)$,
\ie Lison\v{e}k's algorithm. However, using a dedicated implementation 
as in~\cite{handbook} for both point-doubling and point-halving, one would expect Algorithm~\ref{sylow2} to be more
than $150$ times faster than Lison\v{e}k's algorithm (or more than $75$ times
faster with an initial trace check). Since point-doubling for the
dedicated implementation is naturally much faster than MAGMA's, the above time could
be reduced significantly, and Kloosterman zeros for larger fields could
also be found, if required. 

The $O(n)$ factor speedup is due to the fundamental difference between
Lison\v{e}k's algorithm and our approach; while Lison\v{e}k's algorithm
traverses the hypothetically-of-order-$p^n$ Sylow $p$-subgroup tree from leaf to root,
we instead calculate its exact height from root to leaf, which on average is $3$ and 
thus requires an expected single point-halving.

\subsection{Algorithm~\ref{sylow3} for $E_{3^n}(a)$}\label{compare3}

Due to the presence of inversions and square-root computations, one expects
each loop iteration of Algorithm~\ref{sylow3} to be slower than each loop iteration of Algorithm~\ref{sylow2}.
Indeed our basic MAGMA implementation of Algorithm~\ref{sylow3} for curves
defined over $\F_{3^{47}}$ runs $\approx 3.5$ times slower than our one for
Algorithm~\ref{sylow2} for curves defined over $\F_{2^{75}}$. However the
MAGMA implementation is $\approx 15$ times faster than Lison\v{e}k's algorithm
in this case (or equivalently $5$ times faster if a trace check is first
performed).

For the field $\F_{3^{47}} = \F_3[t]/(t^{47} -t^4 - t^2 - t + 1)$, 
using our MAGMA implementation of Algorithm~\ref{sylow3}, we found
the Kloosterman zero:
\begin{eqnarray*}
a & = & t^{46} + t^{45} - t^{44} - t^{42} + t^{39} - t^{38} - t^{36} - t^{35} - t^{33} - t^{31} - t^{30} + t^{29} + t^{28}\\
  & + & t^{26} + t^{25} - t^{24} - t^{22} - t^{21} + t^{20} - t^{19} - t^{17} + t^{16} - t^{15} + t^{14} + t^{13} - t^{11}\\
  & + & t^{10} - t^9 - t^7 + t^6 + t^5 + t^4 -t^2 + 1,
\end{eqnarray*}
in 126 hours, again using eight AMD Opteron 6128 processors running at 2.0 GHz.

In order to improve our basic approach, representational, algorithmic and implementation
optimisations need to be researched. It may be possible for instance to improve the underlying point-thirding algorithm by using
alternative representations of the curve, or the points, or both.
For example, one may instead use the Hessian form~\cite{chudnovsky} of
$E_{3^n}(a)$:
\[
H_{3^n}(\bar{a}): \bar{x}^3 + \bar{y}^3 + 1 = \bar{a}\bar{x}\bar{y},
\]
where $\bar{a} = a^{-1/3}$, $\bar{x} = -a^{1/3}(x+y)$ and $\bar{y} =
a^{1/3}(y-x)$, and an associated tripling formula, see for example~\cite[\S3]{hisil}.
Could point-thirding in this form be faster than that described for the Weierstrass
form in Algorithm~\ref{sylow3}? Also, is there an analogue of the
$\lambda$-representation of a point~\cite{knudsen,schroeppel}
that permits more efficient point-tripling, and hence point-thirding?
We leave as an interesting practical problem the development of efficient 
point-thirding algorithms and implementations for ternary field
elliptic curves with non-zero $j$-invariant.


%

\section{Concluding remarks}\label{conc}

We have presented an efficient deterministic algorithm which tests whether or not
an element of $\F_{2^n}^{\times}$ or $\F_{3^n}^{\times}$ is a Kloosterman zero, and have rigorously analysed
its expected runtime.
Our analysis also gives an upper bound on the number of
Kloosterman zeros. By repeatedly applying our algorithm on random field
elements, we obtain the fastest probabilistic method to date for finding Kloosterman zeros,
which for $\F_{2^n}$ is $O(n)$ times faster than the previous best method, for
$n$ in the practical range. 
Since this method of finding a Kloosterman zero is still exponential in $n$, 
it remains an important open problem to compute Kloosterman zeros efficiently.

\section*{Acknowledgements}

The authors wish to thank Faruk G\"olo\u{g}lu and Alexey Zaytsev for useful
discussions, and the reviewers for their comments.

\bibliographystyle{plain}

\bibliography{KB}

\begin{thebibliography}{10}

\bibitem{omran}
Omran Ahmadi and Alfred Menezes.
\newblock On the number of trace-one elements in polynomial bases for {${\Bbb
  F}_{2^n}$}.
\newblock {\em Des. Codes Cryptogr.}, 37(3):493--507, 2005.

\bibitem{fft}
K.~G. Beauchamp.
\newblock {\em Walsh functions and their applications}.
\newblock Academic Press [Harcourt Brace Jovanovich Publishers], London, 1975.
\newblock Techniques of Physics, No. 3.

\bibitem{ECC1}
I.~F. Blake, G.~Seroussi, and N.~P. Smart.
\newblock {\em Elliptic curves in cryptography}, volume 265 of {\em London
  Mathematical Society Lecture Note Series}.
\newblock Cambridge University Press, Cambridge, 2000.
\newblock Reprint of the 1999 original.

\bibitem{magma}
Wieb Bosma, John Cannon, and Catherine Playoust.
\newblock The {M}agma algebra system. {I}. {T}he user language.
\newblock {\em J. Symbolic Comput.}, 24(3-4):235--265, 1997.
\newblock Computational algebra and number theory (London, 1993).

\bibitem{castryck}
Wouter Castryck and Hendrik Hubrechts.
\newblock The distribution of the number of points modulo an integer on
  elliptic curves over finite fields.
\newblock Preprint, 2011.

\bibitem{gong}
Pascale Charpin and Guang Gong.
\newblock Hyperbent functions, {K}loosterman sums, and {D}ickson polynomials.
\newblock {\em IEEE Trans. Inform. Theory}, 54(9):4230--4238, 2008.

\bibitem{charpin}
Pascale Charpin, Tor Helleseth, and Victor Zinoviev.
\newblock The divisibility modulo 24 of {K}loosterman sums on {${\rm
  GF}(2^m)$}, {$m$} odd.
\newblock {\em J. Combin. Theory Ser. A}, 114(2):322--338, 2007.

\bibitem{charpin2}
Pascale Charpin, Tor Helleseth, and Victor Zinoviev.
\newblock Propagation characteristics of $x \mapsto x^{-1}$ and kloosterman
  sums.
\newblock {\em Finite Fields and Their Applications}, 13(2):366 -- 381, 2007.

\bibitem{chudnovsky}
D.V Chudnovsky and G.V Chudnovsky.
\newblock Sequences of numbers generated by addition in formal groups and new
  primality and factorization tests.
\newblock {\em Advances in Applied Mathematics}, 7(4):385 -- 434, 1986.

\bibitem{dillon}
John.~F. Dillon.
\newblock {\em Elementary Hadamard Difference Sets}.
\newblock PhD Thesis. University of Maryland, 1974.

\bibitem{handbook}
Kenny Fong, Darrel Hankerson, Julio Lopez, and Alfred Menezes.
\newblock Field inversion and point halving revisited.
\newblock {\em IEEE Transactions on Computers}, 53:1047--1059, 2003.

\bibitem{lisonek2}
Kseniya Garaschuk and Petr Lison{\v{e}}k.
\newblock On binary {K}loosterman sums divisible by 3.
\newblock {\em Des. Codes Cryptogr.}, 49(1-3):347--357, 2008.

\bibitem{faruk3}
F.~G{\"o}lo{\u{g}}lu, P.~Lison{\v{e}}k, G.~McGuire, and R.~Moloney.
\newblock Binary kloosterman sums modulo 256 and coefficients of the
  characteristic polynomial.
\newblock {\em Information Theory, IEEE Transactions on}, 58(4):2516--2523,
  2012.

\bibitem{faruk1}
Faruk G{\"o}loglu, Gary McGuire, and Richard Moloney.
\newblock Ternary kloosterman sums modulo 18 using stickelberger's theorem.
\newblock In Claude Carlet and Alexander Pott, editors, {\em SETA}, volume 6338
  of {\em Lecture Notes in Computer Science}, pages 196--203. Springer, 2010.

\bibitem{faruk2}
Faruk G{\"o}lo{\u{g}}lu, Gary McGuire, and Richard Moloney.
\newblock Binary {K}loosterman sums using {S}tickelberger's theorem and the
  {G}ross-{K}oblitz formula.
\newblock {\em Acta Arith.}, 148(3):269--279, 2011.

\bibitem{faruk4}
Faruk G{\"o}loglu, Gary McGuire, and Richard Moloney.
\newblock Ternary kloosterman sums using stickelberger’s theorem and the
  gross-koblitz formula.
\newblock Preprint, 2010.

\bibitem{helleseth1}
Tor Helleseth and Alexander Kholosha.
\newblock Monomial and quadratic bent functions over the finite fields of odd
  characteristic.
\newblock {\em IEEE Transactions on Information Theory}, 52(5):2018--2032,
  2006.

\bibitem{helleseth2}
Tor Helleseth and Victor Zinoviev.
\newblock On {$Z_4$}-linear {G}oethals codes and {K}loosterman sums.
\newblock {\em Des. Codes Cryptogr.}, 17(1-3):269--288, 1999.

\bibitem{hisil}
Huseyin Hisil, Gary Carter, and Ed~Dawson.
\newblock New formulae for efficient elliptic curve arithmetic.
\newblock In K.~Srinathan, C.~Rangan, and Moti Yung, editors, {\em Progress in
  Cryptology – INDOCRYPT 2007}, volume 4859 of {\em Lecture Notes in Computer
  Science}, pages 138--151. Springer Berlin / Heidelberg, 2007.

\bibitem{howe}
Everett~W. Howe.
\newblock On the group orders of elliptic curves over finite fields.
\newblock {\em Compositio Math.}, 85(2):229--247, 1993.

\bibitem{Igusa}
Jun-ichi Igusa.
\newblock On the algebraic theory of elliptic modular functions.
\newblock {\em J. Math. Soc. Japan}, 20:96--106, 1968.

\bibitem{katzM}
N.~M. Katz.
\newblock {\em Gauss sums, Kloosterman sums, and monodromy groups.}
\newblock Princeton Univ. Press, Princeton, NJ, 1988.

\bibitem{katz}
Nicholas Katz and Ron Livn{\'e}.
\newblock Sommes de {K}loosterman et courbes elliptiques universelles en
  caract\'eristiques {$2$} et {$3$}.
\newblock {\em C. R. Acad. Sci. Paris S\'er. I Math.}, 309(11):723--726, 1989.

\bibitem{knudsen}
Erik~Woodward Knudsen.
\newblock Elliptic scalar multiplication using point halving.
\newblock In {\em Proceedings of the International Conference on the Theory and
  Applications of Cryptology and Information Security: Advances in Cryptology},
  ASIACRYPT '99, pages 135--149, London, UK, 1999. Springer-Verlag.

\bibitem{kononen}
K.P. Kononen, M.J. Rinta-aho, and K.O. V\"{a}\"{a}n\"{a}nen.
\newblock On integer values of kloosterman sums.
\newblock {\em Information Theory, IEEE Transactions on}, 56(8):4011 --4013,
  aug. 2010.

\bibitem{wolfmann}
Gilles Lachaud and Jacques Wolfmann.
\newblock The weights of the orthogonals of the extended quadratic binary
  {G}oppa codes.
\newblock {\em IEEE Trans. Inform. Theory}, 36(3):686--692, 1990.

\bibitem{lenstra}
H.~W. Lenstra, Jr.
\newblock Factoring integers with elliptic curves.
\newblock {\em Ann. of Math. (2)}, 126(3):649--673, 1987.

\bibitem{lisonek}
Petr Lison{\v{e}}k.
\newblock On the connection between {K}loosterman sums and elliptic curves.
\newblock In {\em Sequences and their applications---{SETA} 2008}, volume 5203
  of {\em Lecture Notes in Comput. Sci.}, pages 182--187. Springer, Berlin,
  2008.

\bibitem{lisonek3}
Petr Lison{\v{e}}k and Marko Moisio.
\newblock On zeros of kloosterman sums.
\newblock {\em Designs, Codes and Cryptography}, 59:223--230, 2011.
\newblock 10.1007/s10623-010-9457-x.

\bibitem{miret1}
J.~Miret, R.~Moreno, A.~Rio, and M.~Valls.
\newblock Determining the 2-{S}ylow subgroup of an elliptic curve over a finite
  field.
\newblock {\em Math. Comp.}, 74(249):411--427 (electronic), 2005.

\bibitem{miret2}
J.~Miret, R.~Moreno, A.~Rio, and M.~Valls.
\newblock Computing the {$l$}-power torsion of an elliptic curve over a finite
  field.
\newblock {\em Math. Comp.}, 78(267):1767--1786, 2009.

\bibitem{moisio}
Marko Moisio.
\newblock Kloosterman sums, elliptic curves, and irreducible polynomials with
  prescribed trace and norm.
\newblock {\em Acta Arith.}, 132(4):329--350, 2008.

\bibitem{moisio2}
Marko Moisio.
\newblock The divisibility modulo 24 of {K}loosterman sums on {${\rm
  GF}(2^m)$}, {$m$} even.
\newblock {\em Finite Fields Appl.}, 15(2):174--184, 2009.

\bibitem{moisiocode}
Marko Moisio and Kalle Ranto.
\newblock Kloosterman sum identities and low-weight codewords in a cyclic code
  with two zeros.
\newblock {\em Finite Fields Appl.}, 13(4):922--935, 2007.

\bibitem{nied}
Harald Niederreiter.
\newblock The distribution of values of kloosterman sums.
\newblock {\em Archiv der Mathematik}, 56:270--277, 1991.
\newblock 10.1007/BF01190214.

\bibitem{Amilcar}
Am{\'{\i}}lcar Pacheco.
\newblock Rational points on {I}gusa curves and {$L$}-functions of symmetric
  representations.
\newblock {\em J. Number Theory}, 58(2):343--360, 1996.

\bibitem{satoh}
Takakazu Satoh.
\newblock The canonical lift of an ordinary elliptic curve over a finite field
  and its point counting.
\newblock {\em J. Ramanujan Math. Soc.}, 15(4):247--270, 2000.

\bibitem{schroeppel}
R.~Schroeppel.
\newblock Elliptic curves: Twice as fast!
\newblock Presentation at the CRYPTO 2000 Rump Session, 2000.

\bibitem{shpar}
Igor Shparlinski.
\newblock On the values of kloosterman sums.
\newblock {\em IEEE Transactions on Information Theory}, 55(6):2599--2601,
  2009.

\bibitem{Silverman}
Joseph~H. Silverman.
\newblock {\em The arithmetic of elliptic curves}, volume 106 of {\em Graduate
  Texts in Mathematics}.
\newblock Springer-Verlag, New York, 1992.
\newblock Corrected reprint of the 1986 original.

\bibitem{geer}
Gerard van~der Geer and Marcel van~der Vlugt.
\newblock Kloosterman sums and the {$p$}-torsion of certain {J}acobians.
\newblock {\em Math. Ann.}, 290(3):549--563, 1991.

\bibitem{frethesis}
F.~Vercauteren.
\newblock {\em Computing zeta functions of curves over finite fields}.
\newblock PhD Thesis. Katholieke Universiteit Leuven, 2003.

\bibitem{washington}
Lawrence~C. Washington.
\newblock {\em Elliptic curves}.
\newblock Discrete Mathematics and its Applications (Boca Raton). Chapman \&
  Hall/CRC, Boca Raton, FL, 2003.
\newblock Number theory and cryptography.

\end{thebibliography}

\end{document}